\documentclass[11pt,twoside]{article}

\usepackage{latexsym,amssymb,amsthm,amsmath,amscd}

\theoremstyle{plain}

\newtheorem{theorem}{Theorem}[section]

\newtheorem{lemma}{Lemma} [section]

\newtheorem{proposition}{Proposition}[section]

\newtheorem{corollary}{Corollary}[section]

\theoremstyle{remark}
\newtheorem{remark}{Remark} [section]

\theoremstyle{definition}
\def\<{\left < }
\def\>{\right >}
\def\({\left ( }
\def\){\right )}
\def\e{\eqref}
\def\i{\hskip.01in {\rm i}\hskip.01in}
\def\e{\eqref}
\def\p{\partial }

\def \Rt {\tilde{R}}

\newtheorem{example}{Example}[section]
\numberwithin{equation}{section}

\begin{document}
\vbox{\hbox{\small Bull. Transilvania Univ. Brasov  1(50)/2008}
\hbox{\small Proc. Conference RIGA 2008}}
\vskip 1.5truecm

\thispagestyle{empty}

\centerline{\large{\bf On  purely real surfaces in Kaehler surfaces}}

\centerline{\large{\bf  and Lorentz surfaces in Lorentzian Kaehler surfaces}}
\medskip

\centerline{\bf Bang-Yen Chen}

\medskip

 \begin{abstract} {\small An immersion $\phi \colon M \to \tilde M$ of a manifold $M$ into an  indefinite Kaehler  manifold $\tilde M$ is called   purely real if the almost complex structure $J$ on $\tilde M$ carries the tangent bundle of $M$ into a transversal  bundle.
 In this article we survey some recent results on purely real surfaces  in Kaehler surfaces as well as on Lorentz surfaces in Lorentzian Kaehler surfaces. 
 \smallskip
 
}
 \smallskip
 \end{abstract}

\pagestyle{myheadings}

\markboth{B.-Y. Chen}{Purely real and Lorentz surfaces}


\footnote[0]{This article is the written version of my talk delivered at the International Conference on ``Riemannian Geometry and Applications'' held at Bra\c sov, Romania, July 8--11,   2008. The
author would like to express his hearty thanks to the members of the organizing and scientific committees for organizing this  conference dedicated to his 65-th birthday.}

\section{Introduction.}

Let $\tilde M_i^n$ be a complex $n$-dimensional indefinite Kaehler manifold, that means $\tilde M^n_i$ is
endowed with an almost complex structure $J$ and with an indefinite Riemannian
metric $\tilde g$, which is $J$-Hermitian, i.e., for all $p\in \tilde M^n_i$, we have
\begin{align}\label{1.1} &\tilde g(JX, JY)=\tilde g(X,Y), \;\; \forall  X,Y\in T_pM^n,  \\& \label{1.2}\tilde \nabla J=0, \end{align}
where $\tilde \nabla$ is the Levi-Civita connection of $\tilde g$. It follows  that $J$  is
integrable.
The complex index of $\tilde M^n$ is defined as the complex dimension of the largest complex negative
definite subspace of the tangent space. 

When the complex index is one,   the  indefinite Kaehler manifold   $\tilde M^n_1$ is called a {\it Lorentzian Kaehler manifold}.

Let $\tilde M^n_i(4c)$ denote a complete simply-connected indefinite complex space form of constant holomorphic sectional curvature $4c$. The curvature tensor $\tilde R$ of $\tilde M^n_i(4c)$ is given by
\begin{equation}\begin{aligned}\label{1.3}&\tilde R(X,Y)Z=c\{\left<Y,Z\right> \hskip-.02in X\hskip-.02in - \hskip-.02in \left<X,Z\right>\hskip-.02in  Y\hskip-.02in +\hskip-.02in \left<JY,Z\right> JX\\& \ \hskip.7in  -\left<JX,Z\right>\hskip-.02in JY\hskip-.02in +\hskip-.02in 2\left<X,JY\right>\hskip-.02in  JZ\}.\end{aligned}\end{equation}

An immersion $\phi \colon M \to \tilde M$ of a manifold $M$ into an  indefinite Kaehler  manifold $\tilde M$ is called   purely real if the almost complex structure $J$ on the ambient space $\tilde M$ carries the tangent bundle of $M$ into a transversal  bundle, that is  $J(TM)\cap TM=\{0\}$.
  Thus, if  an immersion $\phi \colon M \to \tilde M$ is purely real, it  
  contains no complex points. 
  
  The simplest examples of purely real submanifolds are Lagrangian surfaces and  proper slant  surfaces.

 In this article we survey recent results on purely real surfaces in Kaehler surfaces as well as on Lorentz surfaces in Lorentzian Kaehler surfaces. 
 
 We divide this article into sixteen sections as follows:
 \vskip.01in
 
 \ 1. Introduction.
  \vskip.01in

 \ 2.  Basic formulas and fundamental equations.
  \vskip.01in

\  3. Basics on purely real surfaces.
  \vskip.01in

\  4. Dependence of fundamental equations 
 \vskip.01in

\  5. Purely real minimal surfaces.
 \vskip.01in
 
\  6. A general optimal inequality for purely real surfaces.
  \vskip.01in

\  7. Purely real surfaces in ${\bf C}^2$ satisfying the equality.
 \vskip.01in
 
\  8. Purely real surfaces in $CP^2$ and $CH^2$ satisfying the equality.
  \vskip.01in

\  9. Two general properties for Lorentz surfaces.
  \vskip.01in

  10. Lorentzian complex space forms  and Legendre curves.
 \vskip.01in
 
11. Minimal Lorentz surfaces in ${\bf C}^2_1$.
  \vskip.01in

 12. Two existence results on minimal Lorentz surfaces.
  \vskip.01in

13.  Minimal slant surfaces in $CP^2_1$ and in $CH^2_1$.
  \vskip.01in

  14. Parallel Lorentz surfaces in ${\bf C}^2_1$.
   \vskip.01in

   15. Parallel Lorentz surfaces in $CP^2_1$.
  \vskip.01in
 
  16. Parallel Lorentz surfaces in $CH^2_1$.

\section{Basic formulas and fundamental equations.}

Let $\tilde M^2$ be an  indefinite Kaehler surface. Denote by $\tilde
R$ the Riemann-Christoffel  curvature tensor of $\tilde M^2$. 
Assume  that $M$ is a non-degenerate
surface in $\tilde M^2$ in the sense that the induced metric on
$M$ is non-degenerate. Thus, $M$ is either space-like or Lorentzian.

Denote by $\nabla$ and $\tilde\nabla$ the Levi Civita
connections on $M$ and $\tilde M^2$, respectively. The formulas
of Gauss and Weingarten are given  by (cf.
\cite{c1,O})
\begin{align} &\label{2.1}\tilde \nabla_XY=\nabla_XY+h(X,Y), \\&\label{2.2}
\tilde \nabla_X \xi=-A_\xi X+D_X\xi\end{align} for tangent
vector fields $X,Y$ and a normal vector field $\xi$, where $h$, $A$
and $D$ are the second fundamental form, the shape operator and the
normal connection.

For each $\xi\in T_p^{\perp}M$, the shape operator $A_{\xi}$ is a
symmetric endomorphism of the tangent space $T_pM$ at $p\in M$.

The shape operator and the second fundamental form are related by
\begin{align}\label{2.3} \<h(X,Y),\xi\>=\<A_{\xi}X,Y\>\end{align}
for $X,Y$ tangent to $M$ and $\xi$ normal to $M$. The mean curvature vector  is defined by $H=\frac{1}{2}\, {\rm trace}\,h$.

For a vector $\tilde X\in T_p\tilde M^2$, $p\in M$,  we
denote by $\tilde X^{\top}$ and $\tilde X^\perp$ the
tangential and the normal components of $\tilde X$,
respectively. The equations of Gauss, Codazzi and Ricci are given
respectively by
\begin{align} &(\Rt(X,Y)Z)^{\top} = R(X,Y)Z + A_{h(X,Z)}Y - A_{h(Y,Z)}X, \label{2.4}\\
&(\Rt(X,Y)Z)^{\perp} = (\overline \nabla_X h)(Y,Z)-(\overline\nabla_Y h)(X,Z), \label{2.5}\\ \label{2.6}
&(\Rt(X,Y)\xi)^{\perp} = h(A_{\xi}X,Y)-h(X,A_{\xi}Y)+R^D(X,Y)\xi,\end{align}
for vector fields $X,Y$ and $Z$  tangent to $M$, $\xi$  normal to $M$, where $\overline\nabla h$ and $R^D$ are defined respectively by
\begin{align}\label{2.7}& (\overline\nabla_X h)(Y,Z) = D_X h(Y,Z) - h(\nabla_X Y,Z) - h(Y,\nabla_X Z),\\&\label{2.8} R^D(X,Y)=[D_X,D_Y]-D_{[X,Y]}.\end{align}

The mean curvature vector $\overrightarrow{H}$ and the squared mean curvature $H^2$ of the surface are defined respectively by 
\begin{align}\label{2.9}&\overrightarrow{H}=\frac{1}{2} {\rm trace} \, h,
\\&\label{2.10} H^2=\tilde g(\overrightarrow{H},\overrightarrow{H}).\end{align}

The ellipse of curvature  of a surface $M$ in a Kaehler surface $\tilde M^2$  is
the subset of the normal plane defined as $$\{h(v, v) \in T^\perp_pM:   |v |=1, v\in T_pM,\,p\in M\}.$$
To see that it is an ellipse, we consider an arbitrary orthogonal tangent frame
$\{e_1,e_2\}$. Put $h_{ij}=h(e_i,e_j),i,j=1,2,$ and look at the following formula
\begin{align} \label{2.12}h(v,v)=\overrightarrow{H}+\frac{h_{11}-h_{22}}{2}\cos 2t +h_{12}\sin 2t\end{align}
 for $v=\cos t e_1+\sin t e_2$.
 As $v$ goes once around the unit
tangent circle, $h(v,v)$ goes twice around the ellipse. The ellipse of curvature could degenerate into a line segment or a point. 

 The center of the ellipse is $\overrightarrow H$. The ellipse of curvature is a circle if and only if
the following two conditions hold:\begin{align} \label{2.13} & |h_{11}-h_{22}|^2=4 |h_{12}|^2, \;\;\; \<h_{11}-h_{22},h_{12}\>=0.\end{align}
The property that the ellipse of curvature is a circle is a conformal invariant.

\section{Basics on purely real surfaces.}

An immersion $\phi \colon M \to \tilde M^2$ of a surface $M $ into a Kaehler surface is called {\it purely real} if the almost complex structure $J$ on $\tilde M^2$ carries the tangent bundle of $M$ into a transversal bundle (cf. \cite{c5}, see also \cite{O}). 

A point $p$ on a purely real surface $M$ is called a {\it Lagrangian point} if $J$ carries the tangent space $T_pM$ into its normal space $T^\perp_p M$.
 A purely real surface $M$ is called  {\it  Lagrangian} (or totally real in the sense of \cite{CO}) if the almost complex structure $J$ on $\tilde M^2$ carries the tangent bundle of $M$ into its normal bundle.
  
  For each tangent vector $X$ of a purely real surface $M$, we put
\begin{align} \label{3.1} &JX=PX+FX ,\end{align}
where $PX$ and $FX$ are the tangential and the normal components of $JX$. 

For an oriented orthonormal frame $\{e_1,e_2\}$, it follows from  \e{3.1}  that 
\begin{align} \label{3.2} &Pe_1=(\cos\alpha) e_2,\;\; Pe_2=-(\cos \alpha) e_1\end{align} for some function $\alpha$. This function $\alpha$ is known as the {\it Wirtinger angle}. It is easy to see that the Wirtinger angle is independent of the choice of $e_1,e_2$ which preserves the orientation.

A purely real surface is called a {\it slant surface} if its Wirtinger angle $\alpha$ is constant. For slant surfaces the Wirtinger angle $\alpha$ is called the slant angle  (cf. \cite{c2}).

\section{Dependence of fundamental equations.}

The three fundamental equations of Gauss, Codazzi and Ricci provide the necessary conditions for local isometric embeddability.   These three fundamental equations also play some important roles in physics; 
in particular   in the Kaluza-Klein theory in general relativity theory (cf. \cite{hv,M,S}).

The three fundamental equations of Gauss, Codazzi and Ricci are  independent in general. 
Recently, I am able to prove the following  general property for  purely real surfaces in an arbitrary Kaehler surface \cite{c8}.

\begin{theorem} \label{T:4.1} The equation of Ricci is a consequence of  the equations of Gauss and Codazzi 
for  any purely real surface in any  Kaehler surface.
\end{theorem}

We also have the following.

\begin{theorem} \label{T:4.2}  The equation of Gauss is a consequence of  the equations of Codazzi and Ricci for  purely real surfaces in any Kaehler surfaces.
\end{theorem}

\section{Purely real minimal surfaces in Kaehler surfaces.}

The following result from \cite{c8} provides a necessary condition for  purely real surfaces in complex space forms to be minimal.

\begin{theorem} \label{T:5.1} Let $M$ be a  purely real surface in a complex space form $\tilde M^2(4c)$ of constant holomorphic sectional curvature $4c$. If $M$ is minimal, then the Wirtinger angle $\alpha$ of $M$ satisfies
  \begin{align}\label{5.1} \Delta \alpha=\left\{||\nabla \alpha||^2+6c\sin^2\alpha\right\}\cot\alpha,\end{align}
  where $\Delta$ is the Laplace operator of $M$ and $\nabla\alpha$ is the gradient of $\alpha$.
\end{theorem}

Some easy consequences of Theorem \ref{T:5.1} are the following (see \cite{c8}).

\begin{corollary} \label{C:5.1} \cite{CT} Every   slant surface in  a complex space form $\tilde M^2(4c)$ with $c\ne 0$ is non-minimal unless it is either Lagrangian or complex.
\end{corollary}

\begin{corollary} Every  compact oriented minimal  purely real surface in the complex projective plane $CP^2(4)$ contains some Lagrangian points.
\end{corollary}

\begin{corollary} Let $M$ be a  purely real minimal surface  in ${\bf C}^2$. If the Wirtinger angle $\alpha$ is a harmonic function, then $M$ is slant.
\end{corollary}

\begin{corollary} Let $M$ be a  purely real minimal surface  in $CP^2(4)$. If the Wirtinger angle $\alpha$ is a harmonic function, then $M$ is Lagrangian.
\end{corollary}

A function $f$ on $(M,g)$ is called {\it subharmonic\/} if $\Delta f\geq 0$ holds everywhere on $M$. The surface $M$ is called {\it parabolic\/} if there exists non non-constant negative subharmonic function.

From Theorem \ref{T:5.1} we also have the following.

\begin{corollary} Let $M$ be an oriented minimal  purely real surface in $CP^2(4)$. If $M$ is parabolic, then $M$ contains some Lagrangian points.
\end{corollary}

\section{A general optimal inequality for purely real surfaces.}

For purely real surfaces in complex space forms, we have the following  general optimal inequality proved in \cite{c8}.

\begin{theorem} \label{T:6.1} Let $M$ be a  purely real surface in   a complex space form $\tilde M^2(4c)$. Then we have 
  \begin{align}\label{6.1} H^2\geq 2\{K-||\nabla \alpha||^2-(1+3\cos^2\alpha)c\} +4\<\nabla\alpha, Jh(e_1,e_2)\>\csc \alpha \end{align}
 with respect to any orthonormal frame $\{e_1,e_2\}$ satisfying $\<\nabla \alpha,e_2\>=0$,  where $H^2$ and $K$ denote the squared mean curvature and the Gauss curvature of $M$, respectively.
  
The equality case of \e{6.1} holds at  $p$ if and only if, with respect a suitable adapted orthonormal frame $\{e_1,e_2,e_3,e_4\}$, the shape operators of $M$ at $p$ take the following forms: \begin{align}\label{6.2} A_{e_3}=\begin{pmatrix} 3\varphi & \delta\\ \delta & \varphi \end{pmatrix},\;\;\; A_{e_4}=\begin{pmatrix} \delta+e_1\alpha & \varphi\\ \varphi & 3\delta+3e_1\alpha\end{pmatrix}.\end{align} 
\end{theorem}

When $M$ is a slant surface in $\tilde M^2(4c)$, Theorem \ref{T:6.1} reduces to the following.

\begin{corollary} \label{C:6.1} {\rm  \cite{c4}} If $M$ is a  slant surface in a complex space form $\tilde M^2(4c)$ with slant angle $\theta$, then  we have 
  \begin{align}\label{6.3} H^2\geq 2\{K-(1+3\cos^2\theta)c\} . \end{align}

The equality case of \e{6.3} holds at  $p$ if and only if, with respect a suitable adapted orthonormal frame $\{e_1,e_2,e_3,e_4\}$, the shape operators of $M$ at $p$ take the following forms: \begin{align}\label{16.2} A_{e_3}=\begin{pmatrix} 3\varphi & \delta\\ \delta & \varphi \end{pmatrix},\;\;\; A_{e_4}=\begin{pmatrix} \delta & \varphi\\ \varphi & 3\delta\end{pmatrix}.\end{align} 
 \end{corollary}

The following lemma follows easily from Corollary \ref{C:6.1} and \e{2.12}.

\begin{lemma} \label{L:6.1} If a  slant surface in a complex space form $\tilde M^2(4c)$ satisfies the equality case of  inequality \e{6.3}, then $M$ has circular ellipse of curvature.
\end{lemma}

When the ambient space is ${\bf C}^2$,  Theorem \ref{T:6.1} reduces to the following.

\begin{corollary} \label{C:6.2} Let $M$ be a  purely real surface in   $\bf C^2$. Then  we have 
  \begin{align}\label{6.4} H^2\geq 2\{K-||\nabla \alpha||^2 +2\<\nabla\alpha, J h(e_1,e_2)\>\csc \alpha\} \end{align}
with respect to an orthonormal frame $\{e_1,e_2\}$ satisfying $\<\nabla \alpha,e_2\>=0$.
  
The equality case of \e{6.4} holds if and only if, with respect a suitable adapted orthonormal frame $\{e_1,e_2,e_3,e_4\}$, the shape operators of $M$ take the following forms: \begin{align}\label{6.5} A_{e_3}=\begin{pmatrix} 3\varphi & \delta\\ \delta & \varphi \end{pmatrix},\;\;\; A_{e_4}=\begin{pmatrix} \delta+e_1\alpha & \varphi\\ \varphi & 3\delta+3e_1\alpha\end{pmatrix}.\end{align} 
\end{corollary}

\begin{example} {\rm Let $\alpha(x)$ and $f(y)$ be non-constant real-valued functions and $b$ is a nonzero real number.  Consider the map:
\begin{align}\label{6.6} L(x,y)=\Big(b e^{-\i b^{-1} f(y)}\cot\alpha(x),f(y)\Big).\end{align}
Then the induced metric via \e{6.6} is given by
\begin{align}\notag &g=b^2\alpha'{}^2 \csc^4\alpha dx\otimes dx+ f'{}^2(y)\csc^2\alpha dy\otimes dy,\end{align}
and whose Gauss curvature is  $K=-b^{-2}\sin^4\alpha$. Moreover, it is direct to show that \e{6.6} defines a purely real minimal surface with  Wirtinger angle  $\alpha$ which satisfies the equality case of \e{6.4}.
}\end{example}

 \begin{example} {\rm  Let $\alpha(x)$ and $f(y)$ be non-constant real-valued functions defined on open intervals $I_1,I_2$ respectively and $b$ is a nonzero real number.  
 Consider  $M=(I_1\times I_2,g)$ equipped with metric:
 \begin{align}\label{6.7} & g=b^2\alpha'{}^2(x) \sin^4\alpha(x) dx\otimes dx+f'{}^2(y)\sin^6 \alpha(x) dy\otimes dy. \end{align}
Then the map: 
\begin{align}\label{6.8}\phi(x,y)=\frac{b}{12}\(4 e^{3\i b^{-1}f(y)}\sin^2\alpha(x),\cos (3\alpha(x))-9\cos\alpha(x)\),\end{align}
defines a purely real isometric immersion of $M$ into ${\bf C}^2$ whose Wirtinger angle is $\alpha$.
Moreover, it is direct to show that the squared mean curvature $H^2$,  Gauss curvature $K$, the gradient of  $\alpha$,  and the second fundamental form $h$ of  $M$ satisfy
\begin{equation}\begin{aligned}\label{6.9} &H^2=\frac{4}{b^2\sin^4\alpha},\;\; K=\frac{3}{b^2\sin^4 \alpha},\;\;  \\& ||\nabla\alpha||^2=\frac{1}{b^2\sin^4\alpha},\;\; h\(\frac{\p}{\p x},\frac{\p}{\p y}\)=0.\end{aligned}\end{equation}
Hence, this purely real surface satisfies the equality case of \e{6.4}.
}\end{example}

\begin{example} {\rm Let $w:S^2\to {\bf C}^2$ be the map defined by
$$w(y_0,y_1,y_2)={{1+iy_0}\over {1+y_0^2}}( ry_1,ry_2),\quad y_0^2+y_1^2+y_2^2=1,$$ where $r$ is a positive real number.
Then $w$ is a Lagrangian immersion of the $2$-sphere $S^2$ into ${\bf C}^2$ which is called the {\it Whitney sphere}. 
It is well-known that, up to rigid motions,  the Whitney  sphere is  the only Lagrangian surface in ${\bf C}^2$  satisfying the equality case of \e{6.4} (cf. \cite{cu,c1997}).
}\end{example}

\section{Purely real surfaces in ${\bf C}^2$ satisfying the equality.}

A purely real surface $M$ in a Kaehler surface is said to have {\it full second fundamental form} if its first normal space, ${\rm Im}\, h$, satisfies $\dim\, ({\rm Im}\, h)= 2$ at each point  in $M$. It is said to have
 {\it degenerate second fundamental form} if $\dim\, ({\rm Im}\, h)< 2$ holds at each point in $M$. 

The following result from \cite{c8} classifies purely real minimal  surfaces in ${\bf C}^2$ which satisfy the equality  case of the inequality \e{6.4}.

\begin{theorem}\label{T:7.1} If $M$ is a purely real minimal surface in ${\bf C}^2$  satisfying the equality case of \e{6.4}, then either $M$  is an open portion of a totally geodesic slant plane  or it is congruent to an open portion of 
a negatively curved surface defined by
\begin{align}\label{7.1} L(x,y)=\Big(b e^{-\i b^{-1} f(y)}\cot\alpha(x),f(y)\Big),\end{align}
where $\alpha$ is the Wirtinger angle, $f$ is a non-constant real-valued function, and $b$ is a nonzero real number.  
\end{theorem}

Also, we have following two classification theorems from \cite{c8}.

\begin{theorem}\label{T:7.2} Let $M$ be a purely real surface in ${\bf C}^2$   satisfying the equality case of \e{6.4}. If $M$ has  circular ellipse of curvature, then $M$ is either  an open portion of a  totally geodesic slant plane or an open portion of a Whitney sphere.
\end{theorem}

\begin{theorem}\label{T:7.3} Let $M$ be a purely real surface in ${\bf C}^2$   satisfying the equality case of \e{6.4}. If $M$ has  degenerate second fundamental form, then $M$ is congruent to an open portion of one of the following three types of surfaces$\,:$

\vskip.04in
{\rm (1)} A totally geodesic slant plane.

\vskip.04in
{\rm (2)} A positively curved  surface with Wirtinger angle $\alpha$ defined by
\begin{align}\notag L(x,y)=\frac{b}{12}\(4 e^{3\i b^{-1}f(y)}\sin^2\alpha(x),\cos (3\alpha(x))-9\cos\alpha(x)\),\end{align}
where $\alpha(x)$ and $f(y)$ are non-constant real-valued functions and $b$ is a nonzero real number.  

\vskip.04in
{\rm (3)} A negatively curved surface  with Wirtinger angle $\alpha$ defined by
\begin{align}\notag L(x,y)=\Big(b e^{-\i b^{-1} f(y)}\cot\alpha(x),f(y)\Big),\end{align}
where $\alpha(x)$ and $f(y)$ are non-constant real-valued functions and $b$ is a nonzero real number.  

\end{theorem}

\section{Purely real surfaces in $CP^2$ and $CH^2$ satisfying the equality.}

For purely real surfaces in non-flat complex space form $\tilde M^2(4c)$, $c=1$ or $-1$, we have the following three classification results of B. Y. Chen, A. Mihai and I. Mihai  \cite{cmm}.

\begin{theorem} \label{T:8.1} Let  $M$ be a purely real minimal surface of a complex space form $\tilde M^2(4c)$ with $c=1$ or $-1$. If $M$ satisfies the basic equality, then we have either 
\vskip.05in

{\rm (1)} $c=1$ and $M$ is a totally geodesic Lagrangian surface, or 
\vskip.05in

 {\rm (2)} $c=-1$ and $M$ is a totally geodesic Lagrangian surface  in $CH^2(-4)$, or 
\vskip.05in

{\rm (3)} $c=-1$ and $M$ is congruent to an open portion of a non-slant surface in $CH^2(-4)$ given by the composition $\pi\circ \phi$, where $\pi: H^5_1\to CH^2(-4)$ is the hyperbolic Hopf fibration and
 $\phi:{\bf R}^3\to H^5_1\subset {\bf C}^3_1$ is 
\begin{equation}\begin{aligned} &\notag \phi(x,v,t)=\frac {e^{\i t+\i v/ \sqrt[\uproot{3} 3]{2}}}{3\sqrt[\uproot{3} 3]{4}e^{x}}
\Bigg(3 \i \sqrt[\uproot{3} 3]{2} \sinh \(\tfrac{\sqrt{3}v}{\sqrt[\uproot{3} 3]{2}}\)+\sqrt{3}\big(\sqrt[\uproot{3} 3]{2}+2 e^{2x}\big)\cosh\(\tfrac{\sqrt{3}v}{\sqrt[\uproot{3} 3]{2}}\),  \\&\hskip.1in 3 e^{2x} \cosh\(\tfrac{\sqrt{3}v}{\sqrt[\uproot{3} 3]{2}}\) + \i \sqrt{3}\big(2 \sqrt[\uproot{3} 3]{2}+e^{2x}\big)\sinh \(\tfrac{\sqrt{3}v}{\sqrt[\uproot{3} 3]{2}}\), \frac{\sqrt  {3}}{ e^{3\i v/\sqrt[\uproot{3} 3]{2}}
} \big(\sqrt[\uproot{3} 3]{2}-e^{2x}\big)\Bigg).
\end{aligned}\end{equation}
This purely real  surface has Wirtinger angle $\alpha=\arctan\, (e^{3x})$.
  \end{theorem}

\begin{theorem}\label{T:8.2} Let $M$ be a purely real surface  with circular ellipse of curvature in a complex space form $\tilde M^2(4c)$, $c=1$ or $-1$. If $M$ satisfies the basic equality, then we have either 
\vskip.05in

$(1)$ $M$ is a Lagrangian surface satisfying the equality
\begin{align}\label{8.2} & H^2=2K-2c\end{align}
identically, or
\vskip.05in

$(2)$ $c=-1$ and $M$ is  congruent to an open portion of a proper slant surface in $CH^2(-4)$ given by the composition $\pi\circ \phi$, where $\pi: H^5_1\to CH^2(-4)$ is the Hopf fibration and $\phi: M\to H^5_1\subset {\bf C}^3_1$ is 
\begin{equation}\begin{aligned}&\label{8.3} \phi(u,v,t)=e^{\i t}\(\frac3 2\cosh av+
\frac 1 6{u^2}e^{-av}-\frac \i 6 {{\sqrt{6}u(1+e^{-av})}}-\frac1 2,\right.
\\ &\hskip.2in   \frac1 3(1+2e^{-av})u+\i\sqrt{6}\(-\frac 1 3 +\frac 1 4{e^{av}}+ e^{-av}\(\frac 1
{12}+\frac 1{18}{u^2}\)\), \\ &\hskip.1in   \left. \frac{\sqrt{2}}6 (1-e^{-av})u+\i\sqrt{3}\(\frac 1
6+\frac 1 4{{e^{av}}} +e^{-av}\(-{\frac 5{12}}+\frac 1 {18}{{u^2}}\)\)\)
\end{aligned}\end{equation} with $a=\sqrt{2/ 3}.$
\end{theorem}

\begin{theorem}\label{T:8.3} Let $M$ be a purely real surface   satisfying the basic equality in a complex space form $\tilde M^2(4c)$ with $c=\pm 1$. If $M$ has  degenerate second fundamental form, then either 
\vskip.05in

{\rm (i)}  $M$ is  a totally geodesic Lagrangian surface or 
\vskip.05in

{\rm (ii)} $c=-1$ and $M$ is congruent to an open portion of the surface  in $CH^2(-4)$ given by $\pi\circ \phi$, where $\pi: H^5_1\to CH^2(-4)$ is the hyperbolic Hopf fibration and
 $\phi:{\bf R}^3\to H^5_1\subset {\bf C}^3_1$ is 
\begin{equation}\begin{aligned} &\notag \phi(u,v,t)=\frac {e^{\i t+\i v/ \sqrt[\uproot{3} 3]{2}}}{3\sqrt[\uproot{3} 3]{4}e^{u}}
\Bigg(3 \i \sqrt[\uproot{3} 3]{2} \sinh \(\tfrac{\sqrt{3}v}{\sqrt[\uproot{3} 3]{2}}\)+\sqrt{3}\big(\sqrt[\uproot{3} 3]{2}+2 e^{2u}\big)\cosh\(\tfrac{\sqrt{3}v}{\sqrt[\uproot{3} 3]{2}}\),  \\&\hskip.1in 3 e^{2u} \cosh\(\tfrac{\sqrt{3}v}{\sqrt[\uproot{3} 3]{2}}\) + \i \sqrt{3}\big(2 \sqrt[\uproot{3} 3]{2}+e^{2u}\big)\sinh \(\tfrac{\sqrt{3}v}{\sqrt[\uproot{3} 3]{2}}\),  \frac{\sqrt  {3}}{ e^{3\i v/\sqrt[\uproot{3} 3]{2}}
} \big(\sqrt[\uproot{3} 3]{2}-e^{2u}\big)\Bigg).
\end{aligned}\end{equation}
\end{theorem}

\begin{remark} Lagrangian surfaces in complex space form $\tilde M^2(4c)$ with $c=1$ or $c=-1$ satisfying  equality \e{8.2} have been completely classified in \cite{cu,c1996,CV1}. It follows from Lemma \ref{L:6.1} that such surfaces have circular ellipse of curvature.
\end{remark}

\section{Two general properties for Lorentz surfaces.}

 For arbitrary Lorentz surfaces in arbitrary Lorentzian Kaehler surfaces we have the following two general results.
 
 \begin{theorem} {\rm \cite{c6}} \label{T:9.1} Every Lorentz surface in an arbitrary Lorentzian Kaehler surface is purely real.
 \end{theorem}

\begin{theorem} {\rm \cite{c9}} \label{T:9.2} The equation of Ricci is a consequence of  the equations of Gauss and Codazzi for any  Lorentz surface in any Lorentzian Kaehler surface.
\end{theorem}

\begin{remark}  The converse of Theorem \ref{T:9.1} is false. In fact, not every purely real surface in a Lorentzian Kaehler surface is Lorentzian. A simple such example is the slant plane in ${\bf C}^2_1$  with slant angle $\theta=\cosh^{-1}\sqrt{1+b^2}$ defined by
\begin{align} L(x,y)=\(b y, x+\i \sqrt{1+b^2}y\).\end{align}
where $b$ is a nonzero real number.  It is easy to see that this slant plane is a space-like surface whose  metric tensor  is given by $g=dx^2+dy^2$. 
\end{remark}

\begin{remark} Space-like surfaces in Lorentzian Kaehler surfaces may contain complex points. For instance, $$L(x,y)=(0,x+iy)$$ is a space-like complex plane in ${\bf C}^2_1$.
\end{remark}

\begin{remark}  Theorem \ref{T:9.2} is false in general if the Lorentz surface  were replaced by a spatial surface in a Lorentzian Kaehler surface.
\end{remark}

\section{Lorentzian complex space forms and Legendre curves.}

\subsection{Lorentzian complex space forms} Let {\bf C}$^n$ denote the  complex $n$-plane with complex coordinates
$z_1,\ldots,z_n$. The {\bf C}$^n$ endowed with $g_{i,n}$, i.e., the real part of the Hermitian form
$$b_{i,n}(z,w)=-\sum_{k=1}^i \bar z_kw_k +\sum_{j=i+1}^n \bar z_jw_j,\quad z,w\in\hbox{\bf C}^n,$$ defines a flat indefinite complex space form with complex index $i$. We simply denote the pair ({\bf C}$^n,g_{i,n})$ by {\bf C}$^n_i$. 

Consider the differentiable manifold:
$$S^{2n+1}_2=\{z\in \hbox{\bf C}_1^{n+1} \,;\, b_{1,n+1}(z,z)=1>0 \},$$ which is an indefinite real space form of constant sectional curvature one.   The Hopf fibration
$$\pi: S^{2n+1}_2\to CP^n_1(4):
z\mapsto z\cdot\hbox{\bf C}^*$$ is a submersion and there exists a unique pseudo-Riemannian
metric of complex index one on $CP^n_1(4)$ such that
$\pi$ is a Riemannian submersion. 

The pseudo-Riemannian manifold $CP^n_1(4)$ is a Lorentzian complex space form of constant holomorphic sectional curvature $4$.

Analogously,   consider
$$H^{2n+1}_2=\{z\in\hbox{\bf C}_2^{n+1}\,;\, b_{2,n+1}(z,z)= -1<0 \},$$ which is an
indefinite real space form of constant sectional curvature $-10$. The Hopf fibration
$$\pi: H^{2n+1}_2\to CH^n_1(-4): z\mapsto z\cdot\hbox{\bf C}^*$$ is a submersion
and there exists a unique pseudo-Riemannian metric of complex index 1 on $CH^n_1(-4)$
such that $\pi$ is a Riemannian submersion.

The pseudo-Riemannian manifold $CH^n_1(4c)$ is a Lorentzian complex space form of constant  holomorphic sectional curvature $-4$.

It is well-known that a complete simply-connected Lorentzian complex space form  $\tilde M^n_1(4c)$ is holomorphically isometric to {\bf C}$^n_1$,  $CP^n_1(4)$, or $CH^n_1(-4)$, according to $c=0, \, c=1$ or $c=-1$, respectively (cf. \cite{br}).

\subsection{Special Legendre curves in light cone} 
A  vector $v$ is called {\it space-like} (respectively, {\it time-like}) if  $\<v,v\>>0$ (respectively, $\<v,v\><0$). A vector $v$
is called {\it null} or  {\it light-like} if it is a nonzero vector and it satisfies $\<v,v\>=0$.

The {\it  light cone} $\mathcal LC$ in ${\bf C}^n_i$ ($n\geq 3, i=1,2$) is defined by $$\mathcal LC=\{z\in {\bf C}^n_i:\left<z,z\right>=0\}.$$ A unit speed curve $z(s)$ lying in $\mathcal LC$ is called {\it Legendre} if  $\left<\i z',z\right>=0$ holds identically. For a unit speed Legendre curve $z$ in  $\mathcal LC$, we have $$\left<z,z\right>=\left<z,z'\right>= \left<z,\i z'\right>=\left<\i z,z''\right>= \left<z',z''\right>=0.$$ The Legendre curve $z$  in  $\mathcal LC$ is called  {\it special Legendre} if  $\left<\i z',z''\right>=0$ holds. 
 
The {\it squared curvature} $\kappa^2$ of  a unit speed special Legendre curve $z$ is defined by $$\kappa^2=\left<z'',z''\right>$$ and its {\it Legendre torsion} $\hat \tau$ is defined by $$\hat \tau=\epsilon_z \<z'',\i z'''\>,$$ where $\epsilon_z=1$ or $-1$ according to $z$ is space-like or time-like (cf. \cite{c05,c052}).

\section{Minimal Lorentz surfaces in ${\bf C}^2_1$.}

Lagrangian surfaces in ${\bf C}^2_1$ are Lorentz surfaces automatically. Minimal flat Lagrangian surfaces in the Lorentzian complex plane ${\bf C}^2_1$ have been classified by Chen and L. Vrancken in \cite{CV2}.  On the other hand,   Vrancken proved in \cite{V} that every Lagrangian minimal surface of constant curvature in ${\bf C}^2_1$ is a flat surface (see, also \cite{KV}).

For minimal flat Lorentzian surfaces in the Lorentzian complex plane ${\bf C}^2_1$, we have the following result from \cite{c6}, which can be consideration of an extension of Chen-Vrancken's result in \cite{CV2}.

\begin{theorem}  \label{T:11.1}  Let $\alpha(y)$ and $f(y)$ be two arbitrary differentiable functions of single variable defined on an open interval $I\ni 0$. Then
 \begin{equation}\begin{aligned}\notag & \psi(x,y)=\Bigg(x+\i f(y)+\frac{1}{2} \int_0^y \cosh^2\alpha dy-\int_0^y f'(y)\sinh\alpha dy,\\&\hskip.1in  x-y+\i f(y) +\frac{1}{2} \int_0^y \cosh^2\alpha dy-\int_0^y f'(y)\sinh\alpha dy-\i \int_0^y \sinh\alpha dy\Bigg)\end{aligned}\end{equation}  defines a minimal flat Lorentzian surface in  the Lorentzian complex plane ${\bf C}^2_1$ with $\alpha$ as its Wirtinger angle.

Conversely, every  minimal flat  Lorentzian surface in ${\bf C}^2_1$ is either an open portion of a totally geodesic Lorentzian plane or congruent to the Lorentzian surface described above. 
\end{theorem}

The following result  from \cite{c2008} determines all Lagrangian minimal surfaces in ${\bf C}^2_1$ which are free from flat points.

\begin{theorem}  \label{T:11.2}  Let  $M$ be a  Lagrangian minimal surface in the Lorentzian complex plane ${\bf C}^2_1$. If $M$ contains no flat points, then $M$ is congruent to an open portion of the Lagrangian  surface defined by

\begin{equation}\begin{aligned}\notag & L(x,y)=\Bigg( \int^x_{x_0} \frac{f(x)dx}{\sqrt{f'(x)}}-\frac{\i}{2} \int^x_{x_0} \frac{dx}{\sqrt{f'(x)}} -\i b\int^y_{y_0} \frac{k(y)dy}{\sqrt{k'(y)}}-\frac{ 1}{2b} \int^y_{y_0}  \frac{dy}{\sqrt{k'(y)}}  ,\\&\hskip.4in \int^x_{x_0} \frac{f(x)dx}{\sqrt{f'(x)}}+\frac{\i}{2} \int^x_{x_0} \frac{dx}{\sqrt{f'(x)}}-\i b\int^y_{y_0} \frac{k(y)dy}{\sqrt{k'(y)}}+\frac{1}{2b} \int^y_{y_0}  \frac{dy}{\sqrt{k'(y)}}  \Bigg),
 \end{aligned} \end{equation}
where $b$ is a nonzero real number, $f(x)$ is a  differentiable function with $f'(x)>0$ on an open interval $I_1\ni x_0$, and $k(y)$ is a  differentiable function with $k'(y)>0$ on an open interval $I_2\ni y_0$.
\end{theorem}

The following result from  \cite{c7} completely classifies Lorentzian minimal surfaces in ${\bf C}^2_1$.

\begin{theorem}  \label{T:11.3} Let $z(x)$ and $w(y)$ be two null curves  defined on open intervals $I_1$ and $I_2$ respectively  in the Lorentzian complex plane ${\bf C}^2_1$. If  $\<z(x),w(y)\>\ne 0$ for $(x,y)\in I_1\times I_2$, then 
 \begin{equation}\begin{aligned}\label{11.1} & \psi(x,y)=z(x)+w(y)\end{aligned}\end{equation}  
defines a Lorentzian minimal surface in ${\bf C}^2_1$.

Conversely, locally every Lorentzian minimal surface in ${\bf C}^2_1$ is congruent to a translation surface defined above.
\end{theorem}

It was proved by Chen and J.-M. Morvan  \cite{CM}  that an orientable minimal surface in ${\bf C}^2=({\bf E^4},J)$ is a Lagrangian surface if and only if it is a holomorphic curve with respect to some other orthogonal almost complex structure on ${\bf E}^4$ (for a simple alternate proof of this fact, see \cite{A}).  Using this fact, we know from Theorem \ref{T:11.3} that the situation of Lagrangian minimal surfaces in the Lorentzian complex plane ${\bf C}^2_1$ is quite different from Lagrangian minimal surface in the complex Euclidean plane ${\bf C}^2$.

\section{Two existence results on minimal Lorentz surfaces.}

We have the following two existence results obtained in \cite{c7}.

\begin{proposition} \label{P:12.1} Let $F$ be a nonconstant real-valued function defined on a simply-connected open subset $U$ of ${\bf R}^2$ which satisfies  the following nonlinear Klein-Gordon equation:
\begin{align}\label{12.1} (\ln F)_{uv}=-\frac{1}{F}-F^2.\end{align}
Put $$g_F=- F^{-1}(du\otimes dv+dv\otimes du).$$
  Then, up to rigid motions, there exists a unique Lagrangian minimal immersion $L_F: (U,g_F) \to CP_1^2(4)$ whose second fundamental form satisfies
\begin{equation}\begin{aligned}\notag&  h \text{$\(\frac{\partial}{\partial u},\frac{\partial}{\partial u}\)=F J\frac{\partial}{\partial v}$}, \;\; h\text{ $\(\frac{\partial}{\partial u},\frac{\partial}{\partial v}\)$}=0, \;\;  h\text{ $\(\frac{\partial}{\partial v},\frac{\partial}{\partial v}\)=FJ\frac{\partial}{\partial u}$}.\end{aligned}\end{equation}
\end{proposition}

We call such a Lagrangian minimal surface associated with a solution of the nonlinear Klein-Gordon equation \e{12.1} a {\it Lagrangian minimal  surface of Klein-Gordon type in} $CP^2_1(4)$.

\begin{proposition} \label{P:12.2} Let $P(u,v)$ be a nonconstant real-valued function on a simply-connected open subset $U$ of ${\bf R}^2$ which satisfies  the nonlinear Klein-Gordon equation:
\begin{align}\label{12.3} (\ln P)_{uv}=\frac{1}{P}-P^2.\end{align}
 Put $$g_P=- P^{-1}(du\otimes dv+dv\otimes du).$$
  Then, up to rigid motions, there exists a unique Lagrangian minimal immersion $L_P: (U,g_P) \to CH_1^2(-4)$ whose second fundamental form satisfies
\begin{equation}\begin{aligned}\notag &  h \text{ $\(\frac{\partial}{\partial u},\frac{\partial}{\partial u}\)=P J\frac{\partial}{\partial v}$}, \;\; h\text{ $\(\frac{\partial}{\partial u},\frac{\partial}{\partial v}\)$}=0, \;\;  h\text{ $\(\frac{\partial}{\partial v},\frac{\partial}{\partial v}\)=PJ\frac{\partial}{\partial u}$}.\end{aligned}\end{equation}
\end{proposition}

Similarly, we call a Lagrangian minimal surface in $CH^2_1(-4)$ associated with a solution of the nonlinear Klein-Gordon equation \e{12.3} a {\it Lagrangian minimal surface of Klein-Gordon type in} $CH^2_1(-4)$.

\begin{remark} The two nonlinear Klein-Gordon equations \e{12.1} and \e{12.3} admit infinitely many solutions. Consequently, there are infinitely many Lagrangian minimal surfaces of Klein-Gordon type in $CP^2_1(4)$ and in $CH^2_1(-4)$.
\end{remark}

\section{Minimal slant surfaces in $CP^2_1$ and in $CH^2_1$.}

The following two theorems classify minimal slant surface in the Lorentzian complex projective and complex hyperbolic planes.
 
\begin{theorem}  {\rm  \cite{c7}} \label{T:13.1} Let  $L:M\to CP^2_1(4) $ be a  minimal slant surface in the Lorentzian complex projective plane $CP^2_1(4)$. Then we have:
\vskip.05in

 $(1)$  If $M$ is of constant curvature, then $M$ is congruent to one of the following three types of  surfaces:
\vskip.05in

{\rm (1.a)}  a totally geodesic Lagrangian surface of $CP^2_1(4)$;

\vskip.05in
 {\rm (1.b)} a  curvature one Lagrangian minimal surface defined  by $\pi\circ \tilde L$ with
\begin{equation}\begin{aligned}\notag& \tilde L(x,y)=z'(x)-\frac{2z(x)}{x+y},\end{aligned}\end{equation}
where $z(x), x\in I,$ is a unit speed space-like special Legendre curve  lying in the light cone $\mathcal{LC}\subset {\bf C}^3_1$ with null squared curvature $\kappa^2(s)$, i.e., $\<z''(x),z''(x)\>=0$ on $I$; 

\vskip.05in
 {\rm (1.c)} a flat Lagrangian minimal surface defined by $\pi\circ \tilde L$ with
\begin{equation}\begin{aligned} \notag & \tilde L(x,y)=\text{$ \frac{1}{\sqrt{3}}$}\Bigg(\sqrt{2}e^{\frac{\i }{2a}(x-a^2y)}  \cosh\(\text{$\frac{\sqrt{3}}{2a}$}(x+a^2y)\)  ,e^{\frac{\i}{a} (a^2y-x)}, \\&\hskip1.0in\sqrt{2}e^{\frac{\i }{2a}(x-a^2y)} \sinh\(\text{$\frac{\sqrt{3}}{2a}$}(x+a^2y)\)\Bigg),\end{aligned}\end{equation} where $a$ is a nonzero real number.

\vskip.05in
 $(2)$ If $M$ contains no open subset of constant curvature, then $M$ is  a Lagrangian minimal surface of Klein-Gordon type in $CP^2_1(4)$.
\end{theorem}

\begin{example} {\rm There exist infinitely many unit speed space-like special Legendre curve  lying in the light cone $\mathcal{LC}\subset {\bf C}^3_1$ with null squared curvature. The simplest such examples are the following.
$$z(x)=\(a+\(\frac{1}{4a}+\i b\)s^2,a-\(\frac{1}{4a}-\i b\)s^2, s\),$$
where $a,b$ are nonzero real numbers. It is easy to check that this special Legendre curve has null Legendre torsion, i.e., $\hat \tau=0$.}
\end{example}

\begin{example} {\rm Another example of unit speed space-like special Legendre curve  lying in the light cone $\mathcal{LC}\subset {\bf C}^3_1$ with null squared curvature is the following.
$$\frac{1}{\sqrt{3}}\( e^{\frac{\i }{2}s}\cosh\(\text{\small$\frac{\sqrt{3}s}{2}$}\)- \i \sqrt{3}  e^{\frac{\i}{2}s}\sinh\(\text{\small$\frac{\sqrt{3}s}{2}$}\), 2e^{\frac{\i}{2}s}\sinh\(\text{\small$\frac{\sqrt{3}s}{2}$}\),e^{-\i s}\).$$
This special Legendre curve has constant  Legendre torsion $\hat \tau=-1$.}
\end{example}

\begin{theorem} {\rm \cite{c7}} \label{T:13.2} Let $L:M\to CH^2_1(-4)$ be a  minimal slant surface in the Lorentzian complex hyperbolic plane $CH^2_1(-4)$. Then we have: 
\vskip.05in

 {\rm (1)} If $M$ is of constant curvature, then  $M$ is congruent to one of the following three types of surfaces:
\vskip.05in

{\rm (1.a)}   a totally geodesic Lagrangian surface of $CH^2_1(-4)$;

\vskip.05in
{\rm (1.b)} a Lagrangian minimal surface of constant curvature $-1$ given  by $\pi\circ \tilde L$ with
\begin{equation}\begin{aligned}\notag& \tilde L(x,y)=z'(x)-\sqrt{2}z(x)\tanh\left(\frac{x+y}{\sqrt{2}}\right),\end{aligned}\end{equation}
where $z(x), x\in I,$ is a unit speed time-like special Legendre curve in the light cone $\mathcal{LC}\subset {\bf C}^3_2$ with constant squared curvature $\kappa^2=2$; 

\vskip.05in
{\rm (1.c)} a flat Lagrangian minimal surface defined by $\pi\circ \tilde L$ with
\begin{equation}\begin{aligned} \notag & \tilde L(x,y)=\text{$\frac{1}{\sqrt{3}}$}\Bigg(\sqrt{2}e^{-\frac{\i }{2a}(x+a^2y)}  \cosh\(\text{$\frac{\sqrt{3}}{2a}$}(x-a^2y)\)  ,e^{\i (ay+\frac{x}{a})}, \\&\hskip1.0in\sqrt{2}e^{- \frac{\i }{2a}(x+a^2y)} \sinh\(\text{$\frac{\sqrt{3}}{2a}$}(x-a^2y)\)\Bigg),\end{aligned}\end{equation} where $a$ is a nonzero real number.

\vskip.05in
 $(2)$ If $M$ contains no open subset of constant curvature, then $M$ is  a Lagrangian minimal surface of Klein-Gordon type in $CH^2_1(-4)$.

\end{theorem}

\begin{example} {\rm There exist many unit speed time-like special Legendre curve in the light cone $\mathcal{LC}\subset {\bf C}^3_2$ with constant squared curvature $\kappa^2=2$. The simplest such examples are the following.
$$z(x)=\(\frac{1}{\sqrt{2}},a e^{\sqrt{2}s}-\(\frac{1}{8a}-\i c\)e^{-\sqrt{2}x}, a e^{\sqrt{2}s}+\(\frac{1}{8a}+\i c\)e^{-\sqrt{2}x}\),$$
where $a$ is a nonzero real number.}
\end{example}

 \section{Parallel Lorentz surfaces in ${\bf C}^2_1$.}
 
\subsection{Lorentzian real space forms.} Let ${\bf R}^n_s$ denote the pseudo-Euclidean $n$-space with
metric tensor given by
\begin{align} g_0=-\sum_{i=1}^s dx_i^2+\sum_{j=s+1}^n dx_j^2,\end{align}
where $\{x_1,\ldots,x_n\}$ is the rectangular coordinate system of
$E^n_s$. Then $({\bf R}^n_s,g_0)$ is a flat semi-Riemannian manifold with
index $s$.

We put
\begin{align} & S^n_s=\{x\in {\bf R}^{n+1}_s| \<x,x\>=1\},\\& H^n_s=\{x\in {\bf R}^{n+1}_{s+1}| \<x,x\>=-1\},\end{align}
where $\<\;,\:\>$ is the indefinite inner product on the
pseudo-Euclidean space. It is well-known that $S^n_s$ and $H^n_s$
are complete semi-Riemannian manifolds with index $s$ of constant
sectional curvature 1 and $-1$, respectively. 

The three semi-Riemannian
manifolds ${\bf R}^n_1,S^n_1$ and $H^n_1$ are known as the Minkowski
space-time, the de Sitter space-time and the anti-de Sitter
space-time, respectively. These spaces with index 1 are called
{\it Lorentzian real space forms}.

\subsection{$B$-scroll over null cubic.}

In ${\bf R}^3_1$ take a null curve $\gamma(v)$ with a null frame, i.e., a set  of vector fields $A(v),B(v),C(v)$ such that $$\gamma'(v)=A(v),\;  \<A,A\>=\<B,B\>=\<A,C\>=\<B,C\>=0$$ and $\<A,B\>=\<C,C\>=1$. If these satisfy the following system of equations:
\begin{equation}\begin{aligned} \notag &A'(v)=\kappa_1(v)A(v)+\kappa_2(v)C(v),\; \\&B'(v)=\kappa_1(v)B(v),\; \\&C'(v)=-\kappa_2(v)B(v),\end{aligned}\end{equation}
then $$P(u,v)=\gamma(v)+u B(v)$$ is a Lorentz surface in ${\bf R}^3_1$ which called a $B$-scroll over $\gamma$. 

If $\kappa_1\equiv 0$ and $\kappa_2\equiv 1$ for $\gamma$, then the curve is called the null cubic $C$. In this case, we have the {\it $B$-scroll over the null cubic} $C$ (cf. \cite{magid}).

\subsection{Parallel submanifolds.}
 A submanifold of a pseudo-Riemannian manifold (in particular, in a Riemannian manifold) is called  parallel  if it has
parallel second fundamental form. Parallel submanifolds are one of the most fundamental submanifolds. 

Parallel submanifolds in real  and  complex space forms have been classified in \cite{ferus,T} and in \cite{N1,N2}, respectively. 

\subsection{Classification of parallel Lorentz surfaces.}

Some special classes of parallel submanifolds in Lorentzian space forms  have been studied in \cite{b,g1,g2,magid}.
Complete classification of parallel space-like surfaces and parallel Lorentz surfaces in 4-dimensional Lorentzian real space forms have been obtained in \cite{cv}. 
  
Recently, together with F.  Dillen and J. Van der Veken, we have completely classified parallel Lorentz surfaces in Lorentzian complex space forms. 

The following results are obtained in \cite{cdv}. 

 \begin{proposition}   If $M$ is a parallel Lorentz surface in ${\bf C}^2_1$, then either $M$ is a flat surface or $M$  is totally umbilical in ${\bf C}^2_1$.
 \end{proposition}

 \begin{proposition} If $M$ is a parallel Lorentz surface in $CP^2_1(4)$, then either $M$ is a flat surface or $M$  is a totally geodesic Lagrangian surface.
 \end{proposition}

 \begin{proposition}   If  $M$ is a parallel Lorentz surface in $CH^2_1(4)$, then either $M$ is a flat surface or $M$  is  a totally geodesic Lagrangian surface.
 \end{proposition}

The following result of Chen, Dillen and Van der Veken completely classify parallel Lorentz surfaces in ${\bf R}^4_1$ and also in ${\bf C}^2_1$ (see   \cite{cdv} for details).

 \begin{theorem}   Let $L:M\to {\bf R}^4_2$ be a parallel Lorentz surface in ${\bf R}^4_2$. Then, up to dilations and complex Lorentzian transformations on ${\bf R}^4_2$, $L$ is an open portion of one of the following eleven types of surfaces:
 
 \vskip.05in
$(a)$ A totally geodesic Lorentz plane;

\vskip.05in
$(b)$ The product surface ${\bf R}^1_1\times S^1\subset {\bf R}^2_2\times {\bf R}^2={\bf R}^4_2$, where $S^1\subset {\bf R}^2$ is a circle;

\vskip.05in
$(c)$ The product surface $H^1_1\times {\bf R}^1\subset {\bf R}^2_2\times {\bf R}^2$;

\vskip.05in
$(d)$ The product surface $S^1_1\times {\bf R}^1\subset {\bf R}^2_1\times {\bf R}^2_1$ defined by $$L= (\sinh y,x,0,\cosh y);$$
\vskip.05in

$(e)$ The product surface $H^1\times {\bf R}^1_1\subset {\bf R}^2_1\times {\bf R}^2_1$ defined by $$L= (\cosh x,\sinh x,  y,0);$$

\vskip.05in
$(f)$ A complex circle given by $$L=(c+\i d)(\cos z,\sin z)\in {\bf C}^2_1={\bf R}^4_2$$ with z=x+\i y and  $0\ne c+\i d\in {\bf C}$;

\vskip.05in
$(g)$ A $B$-scroll over a null cubic;

\vskip.05in
$(h)$ A flat surface given by 
\begin{equation}\begin{aligned} &\notag  L=\frac{1}{\sqrt{2}}\Big((1+c)\sin y -(x+cy)\cos y,(1+c)\cos y +(x+cy)\sin y,\\& \hskip.2in(1-c)\sin y +(x+cy)\cos y,(1-c)\cos y -(x+cy)\sin y
\Big),\; c\in {\bf R};
\end{aligned}\end{equation}

\vskip.05in
$(i)$  A marginally trapped flat  surface given by $$L=(q(x,y),x,y,q(x,y)),$$ where $q=ax^2+bxy+cu^2+dx+ey+f$ for some real numbers $a,b,c,d,e,f$ with $a^2+b^2+c^2\ne 0$; 

\vskip.05in
$(j)$ A totally umbilical imbedding of $S^2_1$ in ${\bf R}^4_2$ defined by $$L=\big(0, \sinh x, \cosh x\cos y, \cosh x \sin y\big);$$
\vskip.05in
$(k)$ A totally umbilical imbedding of $H^2_1$ in ${\bf R}^4_2$ defined by $$L=\big(\sin x, \cos x \cosh y, \cos x\sinh y,0\big).$$
\end{theorem}

\begin{remark} By a marginally trapped surface in a Lorentzian Kaehler surface we mean a surface whose mean curvature vector is light-like at each point. Marginally trapped Lagrangian surfaces in Lorentzian complex space forms have been completely classified by Chen and F. Dillen in \cite{cd}. 

Marginally trapped proper slant surfaces in Lorentzian complex space forms were classified by Chen and I. Mihai in \cite{cmihai}. 
Furthermore, marginally trapped  surfaces of constant curvature in Lorentzian complex space forms were completely classified in 
 \cite{c08,c09}.

For a recent survey on marginally trapped surfaces and Kaluza-Klein theory in general relativity, see \cite{mt}.
\end{remark}

\begin{remark} To obtain the classification of parallel Lorentz surfaces in ${\bf C}^2_1$, we only need to equip ${\bf R}^4_2$ with any compatible almost complex structure $J$ on ${\bf R}^4_2$. 
\end{remark}

 \section{Parallel Lorentz surfaces in $CP^2_1$.}

The following result  of Chen, F. Dillen and J. Van der Veken  completely classifies parallel Lorentz surfaces in Lorentzian complex projective plane $CP^2_1$ (see \cite{cdv} for details).

\begin{theorem}  Let $L:M\to CP^2_1(4)$ be a parallel Lorentz surface in $CP^2_1(4)$. Then either 
\vskip.05in

$(a)$ $M$ is an open portion of  the real projective plane $RP^2_1(1)$ of constant curvature one  immersed as a totally geodesic Lagrangian surface in $CP^2_1(4)$, or 
\vskip.05in

$(b)$ $M^2_1$ is a flat surface and up to rigid motions, $L$ is an open portion of the composition $\pi\circ \phi$, where $\pi:S^5_2\to CP^2_1(4)$ is the Hopf fibration and $\phi$ is one of the fourteen types of immersions:

\vskip.05in
$(b.1)$ $\phi:M\to S^5_2\subset {\bf C}^3_1$ is given by
\begin{equation}\begin{aligned} &\notag  \phi=\frac{1}{\sqrt{3}} \(\sqrt{2}e^{\frac{\i}{2}x} \sinh\(\frac{\sqrt{3}}{2}y\), \sqrt{2}e^{\frac{\i}{2}x}\cosh\(\frac{\sqrt{3}}{2}y\),e^{-ix}\);
\end{aligned}\end{equation}

\vskip.05in
$(b.2)$ $\phi:M\to S^5_2\subset {\bf C}^3_1$ is given by
\begin{equation}\begin{aligned} &\notag 
 \phi= \(\frac{e^{\frac{\i}{2} (2x+y+\sqrt{1+4b}y)}} {(1+4b)^{\frac{1}{4} }},  \frac{e^{\frac{\i}{2} (2x+y-\sqrt{1+4b}y)}} {(1+4b)^{\frac{1}{4} }}, e^{\i y}\),\; b>\frac{1}{4};
\end{aligned}\end{equation}

\vskip.05in
$(b.3)$ $\phi:M\to S^5_2\subset {\bf C}^3_1$ is given by
\begin{equation}\begin{aligned} &\notag  \phi=\frac{1}{\sqrt{2}}\(e^{\i (x+\frac{y}{2})}(1+\i y), e^{\i (x+\frac{y}{2})}(1-\i y), \sqrt{2}e^{\i y}\);
\end{aligned}\end{equation}

\vskip.05in
$(b.4)$ $\phi:M\to S^5_2\subset {\bf C}^3_1$ is given by
\begin{equation}\begin{aligned} \notag & \phi=\(\frac{e^{\i (x+\frac{y}{2})}}{(4b-1)^{\frac{1}{4}}}\(\cosh\(\frac{\sqrt{4b-1}}{2}y\)+\i \sinh \(\frac{\sqrt{4b-1}}{2}y\)\),\right.
\\&\hskip-.1in \left. \frac{e^{\i (x+\frac{y}{2})}}{(4b-1)^{\frac{1}{4}}}\(\cosh\(\frac{\sqrt{4b-1}}{2}y\)-\i \sinh \(\frac{\sqrt{4b-1}}{2}y\)\),e^{\i y}\),\; b>\frac{1}{4};
\end{aligned}\end{equation}

\vskip.05in
$(b.5)$ $\phi: M\to S^5_2\subset {\bf C}^3_1$ is given by
\begin{equation}\begin{aligned} &\notag  \phi=\(\frac{e^{\i (bx+\frac{b-1}{2b}y)}}{\sqrt{b-2}} , \frac{e^{\i (-bx+\frac{b-1}{2b}y)}}{\sqrt{b-2}}, \frac{b e^{\i (2x+\frac{y}{b^2})}}{\sqrt{b^2-4}}\),\; 0<b<2;
\end{aligned}\end{equation}

\vskip.05in
$(b.6)$ $\phi:M\to S^5_2\subset {\bf C}^3_1$ is given by
\begin{equation}\begin{aligned} &\notag  \phi=\( \frac{b e^{\i (2x+\frac{y}{b^2})}}{\sqrt{4-b^2}},\frac{e^{\i (bx+\frac{b-1}{2b}y)}}{\sqrt{2-b}} , \frac{e^{\i (-bx+\frac{b-1}{2b}y)}}{\sqrt{2-b}}\),\; b>2;
\end{aligned}\end{equation}

\vskip.05in
$(b.7)$ $\phi:M\to S^5_2\subset {\bf C}^3_1$ is given by
\begin{equation}\begin{aligned}\notag  & \phi=\(   \frac{\sqrt{a(2-a-b)}e^{\i (bx+\frac{(1-b)y}{a(2-a-b)})}}{\sqrt{(a-b)(a+2b-2)}},
   \frac{\sqrt{b(2-a-b)}e^{\i (ax+\frac{(1-a)y}{a(2-a-b)})}}{\sqrt{(a-b)(2a+b-2)}}, \right. \\&\hskip.5in \left. 
 \frac{\sqrt{ab}e^{\i ((2-a-b)x+\frac{a+b-1}{ab}y)}}{\sqrt{(a+2b-2)(2a+b-2)}}  \)
\end{aligned}\end{equation} with $a>b>0,\, a+b<2;$

\vskip.05in
$(b.8)$ $\phi:M\to S^5_2\subset {\bf C}^3_1$ is given by
\begin{equation}\begin{aligned} \notag & \phi=\(   \frac{\sqrt{b(a+b-2)}e^{\i (ax+\frac{(1-a)y}{a(2-a-b)})}}{\sqrt{(a-b)(2a+b-2)}},
    \frac{\sqrt{a(a+b-2)}e^{\i (bx+\frac{(1-b)y}{a(2-a-b)})}}{\sqrt{(a-b)(a+2b-2)}}, \right. \\&\hskip.5in \left. 
 \frac{\sqrt{ab}e^{\i ((2-a-b)x+\frac{a+b-1}{ab}y)}}{\sqrt{(a+2b-2)(2a+b-2)}}  \)
\end{aligned}\end{equation} with $a>b>0,\, a+b>2;$

\vskip.05in
$(b.9)$ $\phi:M\to S^5_2\subset {\bf C}^3_1$ is given by
\begin{equation}\begin{aligned}\notag  & \phi=\(  \frac{\sqrt{-ab}e^{\i ((2-a-b)x+\frac{a+b-1}{ab}y)}}{\sqrt{(a+2b-2)(2a+b-2)}} ,
  \frac{\sqrt{b(2-a-b)}e^{\i (ax+\frac{(1-a)y}{a(2-a-b)})}}{\sqrt{(a-b)(2a+b-2)}},
   \right. \\&\hskip.5in \left.   \frac{\sqrt{a(a+b-2)}e^{\i (bx+\frac{(1-b)y}{a(2-a-b)})}}{\sqrt{(a-b)(a+2b-2)}} \) 
\end{aligned}\end{equation} with $a>2 ,\, b<0;$

\vskip.05in
$(b.10)$ $\phi: M\to S^5_2\subset {\bf C}^3_1$ is given by
\begin{equation}\begin{aligned} &\notag  \phi=\( \(\frac{2\i \sqrt{(2r-1)(1-r)}}{2-3r}+ \frac{2r^2(r-1)x+(2r-1)y}{2r\sqrt{(2r-1)(1-r)}}  \)e^{\i (rx+\frac{y}{2r})} , \right.
\\&\hskip.3in \left.   \frac{2r^2(r-1)x+(2r-1)y}{2r\sqrt{(2r-1)(1-r)}}  e^{\i (rx+\frac{y}{2r})} , \frac{r  e^{\i(2(1-r)x+\frac{2r-1}{r^2}y)}}{3r-2} \)
\end{aligned}\end{equation} with $\frac{2}{3}<r<4;$

\vskip.05in
$(b.11)$ $\phi:M\to S^5_2\subset {\bf C}^3_1$ is given by
\begin{equation}\begin{aligned} &\notag  \phi=\frac{e^{\frac{\i}{12}(8x+9y)}}{288}\Big((8x-9y)^2+24(8x+9y), 24(8x-9y),\\& \hskip.8in 
(8x-9y)^2+24\i (8x+9y+12 \i)\Big);
\end{aligned}\end{equation}

\vskip.05in
$(b.12)$ $\phi: M\to S^5_2\subset {\bf C}^3_1$ is given by
\begin{equation}\begin{aligned} &\notag  \hskip.1in  \phi=\(\frac{e^{-\i rx}}{r(3r-2)}\Big\{2r^2(2+3r^2-5r)x+(2r^2(2-7r+6r^2)t-\i r
\right. \\&\hskip.2in    -\i r^2 \sqrt{5-12r+8r^2}e^{\i(2x+\frac{2r-1}{r^2}y)}\Big\},
 \\&\hskip.0in
\frac{\i r \sqrt{5-12r+8r^2}}{2-3r}e^{\i(2(1-r)x+\frac{2r-1}{r^2}y)}
  +(2r(r-1)x+(2-r^{-1})y)e^{\i( r  x+\frac{y}{2r})}, 
  \\& \left. \hskip.1in  \frac{re^{\i (2(1-r)x+\frac{2r-1}{r^2}y)}}{2-3r}+\frac{\sqrt{5-12r+8r^2}}{2-3r} e^{\i(rx+\frac{y}{2r}y)}  \)
\end{aligned}\end{equation} with $ r\ne 0,1, [\tfrac{2}{3},4);$

\vskip.05in
$(b.13)$ $\phi: M\to S^5_2\subset {\bf C}^3_1$ is given by
\begin{equation}\begin{aligned} \notag &\hskip.0in  \phi=\(\frac{\sqrt{a-1}e^{\i(ax+\frac{(a^2-b^2-a)y}{2(a-1)(a^2+b^2)})}}{b\sqrt{1-2a}\sqrt{(2-3a)^2+b^2}}\big\{ \i(3a^2-b^2-2a)\sinh\(bx+ky\)    \right. \\&\hskip.1in  +(2-4a)b\cosh\(bx+ky\)  \big\},\frac{\sqrt{a^2+b^2}}{\sqrt{(2-3a)^2+b^2}}  e^{\i(2(1-a)x+\frac{2a-1}{a^2+b^2}y)} ,\\ &\left. \hskip.4in\frac{\sqrt{a-1}\sqrt{a^2+b^2}e^{\i(ax+\frac{(a^2-b^2-a)y}{2(a-1)(a^2+b^2)})}}{b\sqrt{1-2a}}\sinh(bx+ky)\)
\end{aligned}\end{equation}
with $$k=\frac{(2a-1)b}{2(a-1)(a^2+b^2)},\;\; a\in ( \tfrac{1}{2},1);$$

\vskip.05in
$(b.14)$ $\phi: M\to S^5_2\subset {\bf C}^3_1$ is given by
\begin{equation}\begin{aligned}\notag  &\hskip.0in  \phi=\(\frac{\sqrt{a-1}e^{\i(ax+\frac{(a^2-b^2-a)y}{2(a-1)(a^2+b^2)})}}{b\sqrt{2a-1}\sqrt{(2-3a)^2+b^2}}\big\{ \i(3a^2-b^2-2a)\cosh\(bx+ky\)    \right.
\\&\hskip.1in  +(2-4a)b\sinh\(bx+ky\)  \big\},\frac{\sqrt{a^2+b^2}}{\sqrt{(2-3a)^2+b^2}}  e^{\i(2(1-a)x+\frac{2a-1}{a^2+b^2}y)} ,\\&\left. \hskip.4in\frac{\sqrt{a-1}\sqrt{a^2+b^2}e^{\i(ax+\frac{(a^2-b^2-a)y}{2(a-1)(a^2+b^2)})}}{b\sqrt{2a-1}}\cosh(bx+ky)\)
\end{aligned}\end{equation}
with $$k=\frac{(2a-1)b}{2(a-1)(a^2+b^2)},\;\; a\notin [ \tfrac{1}{2},1].$$

 \end{theorem}

 \section{Parallel Lorentz surfaces in $CH^2_1$.}

The following result  of Chen, F. Dillen and J. Van der Veken  completely classifies parallel Lorentz surfaces in Lorentzian complex hyperbolic plane $CH^2_1$ (see \cite{cdv} for details).

\begin{theorem}  Let $L:M\to CH^2_1(-4)$ be a parallel Lorentz surface in $CH^2_1(-4)$. Then either 
\vskip.05in

$(i)$ $M$ is an open portion of the real hyperbolic plane $H^2_1(-1)$ of constant curvature $-1$ immersed as a totally geodesic Lagrangian surface in $CH^2_1(-4)$, or 
\vskip.05in

$(ii)$ $M^2_1$ is a flat surface and up to rigid motions, $L$ is an open portion of the composition $\pi\circ \phi$, where $\pi:H^5_2\subset  CH^2_1(-4)$ is the Hopf fibration and $\phi$ is one of the thirteen types of immersions:

\vskip.05in
$(ii.1)$ $\phi:M\to H^5_2\subset  {\bf C}^3_2$ is given by
$\phi=\(e^{\i y},e^{\i x},e^{\i(x+y)}\);$

\vskip.05in
$(ii.2)$ $\phi:M\to H^5_2\subset {\bf C}^3_2$ is given by
\begin{equation}\begin{aligned}\notag  & \phi=\(\(y-\frac{\i}{2}\)e^{\i(x-\frac{y}{2})}, e^{-\i y}, \(y+\frac{\i}{2}\)e^{\i(x-\frac{y}{2})}\);
\end{aligned}\end{equation}

\vskip.05in
$(ii.3)$ $\phi:M\to H^5_2\subset  {\bf C}^3_2$ is given by
\begin{equation}\begin{aligned} &\notag  \phi=\frac{1}{\sqrt{6}} \(\frac{e^{\sqrt{3}x+\sqrt{3}y}}{e^{\frac{1}{2}(\i+\sqrt{3})(x+y)}},
\frac{e^{\sqrt{3}x-\sqrt{3}y}}{e^{\frac{1}{2}(\i+\sqrt{3})(x+y)}}, \i \sqrt{2}e^{\i(x+y)}\);
\end{aligned}\end{equation}

\vskip.05in
$(ii.4)$ $\phi:M\to H^5_2\subset  {\bf C}^3_2$ is given by
\begin{equation}\begin{aligned} &\notag  \phi=\( e^{-\i y}, \frac{e^{\i (x+\frac{y}{2}(\sqrt{1+4b}-1)}}{(1+4b)^{\frac{1}{4}}} ,  \frac{e^{\i (x-\frac{y}{2}(\sqrt{1+4b}+1)}}{(1+4b)^{\frac{1}{4}}}  \),\; b>-\frac{1}{4};
\end{aligned}\end{equation} 

\vskip.05in
$(ii.5)$ $\phi:M\to H^5_2\subset  {\bf C}^3_2$ is given by
\begin{equation}\begin{aligned} &\notag\;  \phi=\( e^{-\i y}, \frac{(2 e^{\sqrt{4b-1}y}-\i)e^{\i x-\frac{y}{2}(\i +\sqrt{4b-1})}}{2(4b-1)^{\frac{1}{4}}} ,  \frac{(2 e^{\sqrt{4b-1}y}+\i)e^{\i x-\frac{y}{2}(\i+\sqrt{4b-1})}}{2(4b-1)^{\frac{1}{4}}}  \)
\end{aligned}\end{equation}
with $b>-\frac{1}{4};$
\vskip.05in
$(ii.6)$ $\phi:M\to H^5_2\subset  {\bf C}^3_2$ is given by
\begin{equation}\begin{aligned}\notag  & \phi=\(   \frac{\sqrt{b(2-a-b)}e^{\i (ax+\frac{(a-1)y}{b(2-a-b)})}}{\sqrt{(a-b)(2a+b-2)}},  \frac{\sqrt{ab}e^{\i ((2-a-b)x+\frac{1-a-b}{ab}y)}}{\sqrt{(a+2b-2)(2a+b-2)}}, \right. \\&\hskip.2in \left.   \frac{\sqrt{a(2-a-b)}e^{\i (bx+\frac{(b-1)y}{a(2-a-b)})}}{\sqrt{(a-b)(a+2b-2)}}  \),\; a>b>0,\, a+b<2;
\end{aligned}\end{equation} 

\vskip.05in
$(ii.7)$ $\phi:M\to H^5_2\subset  {\bf C}^3_2$ is given by
\begin{equation}\begin{aligned}\notag  & \phi=\(   \frac{\sqrt{a(a+b-2)}e^{\i (bx+\frac{(b-1)y}{a(2-a-b)})}}{\sqrt{(a-b)(2a+b-2)}}, \frac{\sqrt{ab}e^{\i ((2-a-b)x+\frac{1-a-b}{ab}y)}}{\sqrt{(a+2b-2)(2a+b-2)}} 
  , \right. \\&\hskip.5in \left.  \frac{\sqrt{b(2-a-b)}e^{\i (ax+\frac{(a-1)y}{b(2-a-b)})}}{\sqrt{(a-b)(2-2a-b)}} \),\; a>b>0,\; a+b>2;
\end{aligned}\end{equation}

\vskip.05in
$(ii.8)$ $\phi:M\to H^5_2\subset  {\bf C}^3_2$ is given by
\begin{equation}\begin{aligned}\notag  & \phi=\(    \frac{\sqrt{b(2-a-b)}e^{\i (ax+\frac{(a-1)y}{b(2-a-b)})}}{\sqrt{(a-b)(2a+b-2)}}, 
  \frac{\sqrt{a(a+b-2)}e^{\i (bx+\frac{(b-1)y}{a(2-a-b)})}}{\sqrt{(a-b)(a+2b-2)}},
   \right. \\&\hskip.5in \left.   \frac{\sqrt{-ab}e^{\i ((2-a-b)x+\frac{1-a-b}{ab}y)}}{\sqrt{(a+2b-2)(2a+b-2)}} \) 
\end{aligned}\end{equation} with $a>0 ,\,b<0,\, a+2b>2;$

\vskip.05in
$(ii.9)$ $\phi:M\to H^5_2\subset  {\bf C}^3_2$ is given by
\begin{equation}\begin{aligned} &\notag  \phi=\(  \frac{r  e^{\i(2(1-r)x+\frac{1-2r}{r^2}y)}}{3r-2}, \right.
  \frac{2r^2(r-1)x+(1-2r)y}{2r\sqrt{(2r-1)(1-r)}}  e^{\i (rx-\frac{y}{2r})} ,\\&\hskip.0in \left.  \(\frac{2\i \sqrt{(2r-1)(1-r)}}{2-3r}+ \frac{2r^2(r-1)x+(1-2r)y}{2r\sqrt{(2r-1)(1-r)}}  \)e^{\i (rx-\frac{y}{2r})}\)
\end{aligned}\end{equation}
with $r\in (\frac{1}{2},1)$ and $r\ne \frac{2}{3};$
\vskip.05in

$(ii.10)$ $\phi:M\to H^5_2\subset  {\bf C}^3_2$ is given by
\begin{equation}\begin{aligned} &\notag  \hskip.1in  \phi= \( \frac{ 2r(r-1)\{2\i (2r-1)+r(3r-2)x+(7r-2-6r^2)y\}}{2r(2-3r)\sqrt{1-3r+2r^2}}e^{\i (rx-\frac{y}{2r} )}  , 
\right. \\&\hskip.3in \left. \frac{r}{3r-2} e^{\i(2(1-r)x+\frac{1-2r}{r^2}y)}, \frac{ 2r^2(r-1) x+2r(2r-1)y}{2r\sqrt{1-3r+2r^2}}e^{\i (rx-\frac{y}{2r} )}  \)\end{aligned}\end{equation}
with $r\notin[\frac{1}{2},1]$ and $r\ne 2$;

\vskip.05in
$(ii.11)$ $\phi:M\to H^5_2\subset  {\bf C}^3_2$ is given by
\begin{equation}\begin{aligned} &\notag  \phi=\frac{e^{\frac{\i}{12}(8x-9y)}}{24}\Big((8x+9y)^2-432\i y-1, 16x+18y+24\i,\\& \hskip.8in  (8x-9y)^2-431\i y+1)\Big);\end{aligned}\end{equation}

\vskip.05in
$(ii.12)$ $\phi:M\to H^5_2\subset  {\bf C}^3_2$ is given by
\begin{equation}\begin{aligned} \notag &\hskip.0in  \phi=\(\frac{\sqrt{a-1}e^{\i(ax+\frac{(a-a^2+b^2)y}{2(a-1)(a^2+b^2)})}}{b\sqrt{2a-1}\sqrt{(2-3a)^2+b^2}}\big\{ \i(2a-3a^2+b^2)\sinh\(bx+ky\)    \right.
\\&\hskip.1in  +(4a-2)b\cosh\(bx+ky\)  \big\},
\frac{\sqrt{a^2+b^2}}{\sqrt{(2-3a)^2+b^2}}  e^{\i(2(1-a)x+\frac{1-2a}{a^2+b^2}y)} ,
\\&\left. \hskip.6in\frac{\sqrt{a-1}\sqrt{a^2+b^2}}{b\sqrt{2a-1}}e^{\i(ax+\frac{(a-a^2+b^2)y}{2(a-1)(a^2+b^2)})}\sinh(bx+ky)\)
\end{aligned}\end{equation}
with $$k=\frac{(1-2a)b}{2(a-1)(a^2+b^2)},\;\; a\notin [ \tfrac{1}{2},1);$$

\vskip.05in
$(ii.13)$ $\phi:M\to H^5_2\subset  {\bf C}^3_2$ is given by 
\begin{equation}\begin{aligned} \notag &\hskip.0in  \phi=\(\frac{\sqrt{a^2+b^2}}{\sqrt{(2-3a)^2+b^2}}  e^{\i(2(1-a)x+\frac{1-2a}{a^2+b^2}y)},\right. \\&\hskip.3in\frac{\sqrt{1-a}e^{\i(ax+\frac{(a-a^2+b^2)y}{2(a-1)(a^2+b^2)})}}{b\sqrt{2a-1}\sqrt{(2-3a)^2+b^2}}\Big\{ (4a-2)b\sinh\(bx+ky\) 
\\&\hskip.7in  + \i(2a-3a^2+b^2)\cosh\(bx+ky\)   \Big\},
\\&\left. \hskip.4in\frac{\sqrt{a-1}\sqrt{a^2+b^2}}{b\sqrt{2a-1}}e^{\i(ax+\frac{(a-a^2+b^2)y}{2(a-1)(a^2+b^2)})}\cosh(bx+ky)\)\end{aligned}\end{equation}
with $$k=\frac{(2a-1)b}{2(a-1)(a^2+b^2)},\;\; a\in ( \tfrac{1}{2},1).$$
 \end{theorem}

\vskip.2in

\noindent {\small Department of Mathematics,}

\noindent {\small  Michigan State University, }
 
\noindent {\small  East Lansing, Michigan 48824-1027,}

\noindent {\small  U.S.A. }

\noindent {\small {\it Email address:} bychen@math.msu.edu}
\end{document}